\title{The Candy-Passing Game for $c\geq3n-2$}

\author{Paul M. Kominers\thanks{Walt Whitman High School, Bethesda, MD}}
\documentclass[12pt]{article}
\usepackage{amsmath}
\usepackage{amsthm}
\usepackage{amsfonts}
\usepackage{mathenv}
\usepackage{setspace}
\usepackage{mathenv}
\newtheorem{thm}{Theorem}
\newtheorem{lemma}{Lemma}

\theoremstyle{remark}

\begin{document}
\maketitle
\begin{abstract}
We determine the behavior of Tanton's candy-passing game for all distributions of at least $3n-2$ candies, where $n$ is the number of students.  Specifically, we show that the configuration of candy in such a game eventually becomes fixed.
\end{abstract}
The candy-passing game, as introduced by Tanton \cite{Candy Passing}, is played according to the following rules:

\begin{itemize}
\item At the beginning of the game, $c>0$ candies are distributed
  arbitrarily among $n>2$ students, who are sitting in a circle.
\item A whistle is sounded at a regular interval.
\item Each time the whistle is sounded, each student with two or more
  candies passes one candy to his left-hand neighbor and one candy to his right-hand neighbor.
\item If a student has fewer than two candies when the whistle sounds, he
  does nothing.
\end{itemize}

When $c<n$, the game eventually terminates, with no students having sufficient candy to pass candy to their neighbors (see \cite{Candy Passing}).  In this paper, we study the behavior of the candy-passing game for $c\geq3n-2$, showing that the configuration of candy in such a game eventually becomes fixed.

We call the interval between blows of the whistle a \emph{round} of candy-passing.  The students are consecutively numbered $1, 2,\ldots,k,\ldots,n$, starting from an arbitrary student.  The candy pile of a student having four or more pieces of candy is said to be \emph{abundant}, and we denote the number of students with abundant candy piles by $m$.  If, after some round, the amount of candy a student has will not change over the remainder of the given candy-passing game, that student's candy pile is said to have \emph{stabilized}.

Clearly, if a student has two or more candies at the beginning of a round, that student cannot end the round with more candy than he began with.  Indeed, in any round, a given student with two or more candies can, at most, pass two pieces of candy to his neighbors and receive two pieces of candy from his neighbors, resulting in no net increase in the size of his candy pile.

\begin{lemma}
\label{Fixed Depletions}
After a finite number of rounds, the set of students with abundant candy piles in any candy-passing game is fixed and the candy piles of all such students have stabilized.
\end{lemma}

\begin{proof}
If $m=0$ at the beginning of the game, there are no students with abundant candy piles to lose candy.  We now assume that $m>0$ at the beginning of the game.  As we observed, the total amount of candy possessed by students with abundant candy piles is nonincreasing.  Further, if a student with an abundant candy pile loses candy, that sum decreases.  Since the total amount of candy possessed by students with abundant candy piles cannot fall below zero, the amount of candy that can be lost by students with abundant candy piles must be finite.
\end{proof}

We are now ready to prove our main result.

\begin{thm}
\label{Annoying result}
In a candy-passing game with $c\geq{3n-2}$, then all students' candy piles eventually stabilize.
\end{thm}

\begin{proof}
As a consequence of Lemma \ref{Fixed Depletions}, we may assume that all candy that may be lost by students with abundant candy piles has been lost, as this must happen within finitely
many rounds.  If $m=0$ at this point, then the condition $c\geq{3n-2}$ implies $c=3n$, $c=3n-1$, or $c=3n-2$.

If $m=0$ and $c=3n$, each student has three candies.  If $m=0$ and $c=3n-1$, each student has three candies except for one, who has two.  In both of these cases, all of the students' candy piles have stabilized.

If $m=0$ and $c=3n-2$, either each student has three candies except for two students who have two candies each, or each student has three candies except for one student who has only one candy.  In the first of these cases, all of the students' candy piles have stabilized.  In the second, the neighbors of the student having only one candy each pass him one candy and receive one, reducing this situation to the first case.

We now assume $m>0$.  Since a student with an abundant candy pile passes candy each round, in order for his candy pile to have stabilized, he must be receiving candy from both of his neighbors every round.  Select one student in the game with an abundant candy pile.  Since he must pass candy every round, he must receive candy from both of his neighbors every round, who must therefore themselves have at least two pieces of candy every round.  The neighbors must therefore eventually stabilize; there is a minimum amount of candy (two pieces) such that they can pass candy every round.  They cannot drop below this number, or the abundant candy piles would not have stabilized.  They cannot oscillate between various amounts of candy greater than two (say, between two and three), as any student with two or more pieces of candy cannot end the round with more candy than he began with.  For them to have stabilized while passing candy every round, their neighbors must be passing candy every round, which means that they, too, must eventually stabilize.  We see by this argument that for an abundant candy pile to have stabilized, all other candy piles must eventually stabilize.
\end{proof}

\paragraph{\textbf{Acknowledgments.}}
This research was conducted at the 2007 Research Science Institute at the Massachusetts Institute of Technology and was supported by a fellowship from the Center for Excellence in Education.

The author is especially grateful to Amanda Redlich and Allison Gilmore, who supervised his research, and to his math teacher, Susan Schwartz Wildstrom, for her encouragement.

He also appreciates the helpful comments and suggestions of Zachary Abel, Gabriel D. Carroll, Noam D. Elkies, Scott D. Kominers, James Propp, John Rickert, Robert Obryk, Jenny Sendova, Peter Shor, James Tanton, and an anonymous referee.


\begin{thebibliography}{99}

\bibitem{Candy Passing} J. Tanton.  Today's puzzler.  \emph{The St. Mark's Institute of Mathematics Newsletter}.  November 2006.

\end{thebibliography}
\end{document}